\documentclass[proc]{edpsmath}
\usepackage{dsfont}
\usepackage{graphicx}
\usepackage{hyperref}

\def\norm#1{\left\| \mskip 1mu #1 \mskip 1mu \right\|}

\def\vgx{\mathbf{x}}
\def\vgy{\mathbf{y}}
\def\vgf{\mathbf{f}}
\def\vga{\mathbf{a}}
\def\vgt{\mathbf{t}}
\def\vgp{\mathbf{p}}
\def\vgK{\mathbf{K}}
\def\vgM{\mathbf{M}}
\def\vgR{\mathbf{R}}
\def\vgS{\mathbf{S}}
\def\vgQ{\mathbf{Q}}
\def\vgP{\mathbf{P}}
\def\egg{\boldsymbol{\varepsilon}}
\def\ngg{\boldsymbol{\nu}}
\def\e{\varepsilon}

\begin{document}

\title{Iterative reconstruction methods for wave equations}
\author{S\'ebastien Marinesque}\address{Universit\'e Paul Sabatier, Institut de Math\'ematiques de Toulouse, 118 route de Narbonne, F-31062 Toulouse Cedex}
%
%
\begin{abstract} 
Some iterative techniques are defined to solve reversible inverse problems and a common formulation is explained. Numerical improvements are suggested and tests validate the methods.
\end{abstract}

\begin{resume}
Nous d\'efinissons des techniques it\'eratives pour la r\'esolution de probl\`emes inverses r\'eversibles et en fournissons une formulation commune. Apr\`es avoir sugg\'er\'e des am\'eliorations pour leur impl\'ementation, des exp\'erimentations sont pr\'esent\'ees qui valident ces m\'ethodes. 
\end{resume}
\maketitle
\section*{Introduction}
The classical ThermoAcoustic Tomography (TAT) problem formulates as follows. Consider an object~$f_0$ contained in an open set~$\Omega$ of~$\xR^n$ which emits an acoustic pressure wave at~$t=0$, considered as a Dirac pulse. This wave is modeled as a solution of:
\begin{equation}
\label{eqn}
\left\{
\begin{array}{l}
L p=0\text{ in }(0,T)\times\Omega,\\
p(0)=f_0,\\
\partial_t p(0)=0,
\end{array}
\right.
\end{equation}
where~$L$ is an operator modeling an acoustic wave phenomenon. Then this pressure wave is observed (e.g. thanks to piezoelectric sensors) and a set of observations is obtained from the solution~$p$. That can be expressed thanks to an observation operator~$C$ mapping a solution~$p$ to observations~$Cp$. The inverse TAT problem consists in developing and studying methods to reconstruct~$f_0$ from~$Cp$ and to define situations in which this reconstruction is possible.

In the three past decades many techniques have been developed, offering effectual results (see works from authors of~\cite{akuku},~\cite{finch}~\cite{aust},~\cite{qsuz},~\cite{xw} among others). The new techniques we propose in section~\ref{iswe} rely on the following idea, influenced by~\cite{AB}: if the system we consider is reversible in time, then the initial state to reconstruct can be seen, backward in time, as a state to reach, so that usual control and filtering techniques can be used to solve this inverse problem. For this purpose, we first used the Back and Forth Nudging algorithm (see~\cite{AB}) in~\cite{oun}. With filtering techniques, as the Kalman filter defined in~\cite{kal} and one of its reduced rank formulation, the SEEK filter (cf.~\cite{pham}), we deal here with possible improvements of this method.

Of course, many assumptions are necessary to obtain a favorable observation situation: the way the wave propagates, depending on the media and on the kind of the wave, the final time and the number, size and position of the sensors recording the wave command the information contained in the data (see e.g.~\cite{kuku},~\cite{liu} and the references therein). Moreover, even if the continuous problem is well set, numerical issues still put up some resistance, as considering noisy data, algorithmic complexity or apparition of spurious high-frequency oscillations during numerical implementations (the paper~\cite{zua05} surveys this point).

We deal here with this issue, as done in~\cite{rtt}. We introduce an artificial attenuation term in the numerical scheme that not only yields a regularization of the solution but corrects degenerated observation configurations for which filtering techniques are not helpful.

\section{Iterative stabilization and filtering for wave equations}\label{iswe}
This section is devoted to main results about stabilization of the wave equation and Kalman-Bucy filter. Then we define iterative stabilizing methods.
\subsection{Observation and stabilization of the wave equation}\label{osw}
We assume that~$L=\partial_{tt}-\Delta$ is the D'Alembert operator. If necessary, we can suppose that~$f_0=0$ and consider the stabilization problem for the initial value problem related to~$L$ instead of the inverse TAT problem. Similar considerations are also valid for linear variable speed (reversible) wave equations.

Define~$\vgp=\left(\begin{array}{c}p\\p'\end{array}\right)$, for any~$p\in\mathcal{C}^0((0,T);H^1_0)\cap\mathcal{C}^1((0,T);L^2)$, and~$A=\left(\begin{array}{cc}0&I\\\Delta&0\end{array}\right)$ on~$H=H^1_0\times L^2$, then~$D(A)=(H^1_0\cap H^2)\times H^1_0$ and equation~\eqref{eqn} writes:
\begin{equation}
\label{eqn1}
\left\{
\begin{array}{l}
\vgp'=A\vgp,\\
\vgp(0)=\vgp_0=\left(\begin{array}{c} f_0\\0\end{array}\right).
\end{array}
\right.
\end{equation}


The observations are defined in a Hilbert space~$U$. It is convenient, when working with wave equations, to consider the time derivative of the observations, so that we assume that~$C\in\mathcal{L}(L^2,U)$ and use equally~$Cp'$ or~$Cp$. In the practice of TAT, we only get observations from~$p$ and use its time derivative when needed.

Concerning stabilizability and controllability of wave equations, we have the fundamental criterion:  
\begin{dfntn} The \emph{observation inequality} is satisfied if there exists~$T,M > 0$ such that:
\begin{equation}
\label{oc}
\int_0^T \norm{C p'}_U^2 dt\geq M\norm{\vgp_0}^2_H,
\end{equation}
for all~$\vgp_0\in H$, where~$p$ is the solution of~\eqref{eqn} with initial data~$\vgp_0$. Scalar~$M$ is called \emph{observability constant}.
\end{dfntn}

Indeed, in~\cite{liu}, one finds the following result:
\begin{prpstn}\label{prop1}The three following propositions are equivalent:
\begin{itemize}
\item[($i$)] The observation inequality~\eqref{oc} is satisfied.
\item[($ii$)] For every positive-definite self-adjoint operator~$T\in L(U)$, the operator~$A-C^\star TC$ generates an exponentially stable~$\mathrm{C}_0$-semigroup on~$H$.
\item[($iii$)] The system~$(A,C^\star)$ is exactly controllable. 
\end{itemize}
\end{prpstn}

Many geometrical interpretations to the observation inequality have been presented, mostly known as \emph{Geometric Optics Condition} (GOC from~\cite{blr}). In particular, these results explain this heuristic situation: when~$Cp'=\mathds{1}^{\vphantom{+}}_\omega p'$, where~$\omega$ is an open subset of~$\Omega$, then enough energy from~$\vgp_0$ has to pass through~$\omega$ to get enough information to reconstruct~$\vgp_0$. It depends on many parameters such as the position of the sensors, the speed map of the wave equation, the final time~$T$, etc.

In this context, we introduce a first reconstruction method, the Kalman-Bucy filter.

\subsection{The Kalman-Bucy filter}\label{kseek}
We recall the main results concerning the (continuous) Kalman-Bucy filter. It yields a way to approximate the real state that minimizes the error variance in the following situation (see~\cite{kb}):

Assume that the true state solves the linear differential equation~${\vgx^\vgt}' = \vgM\vgx^\vgt+\ngg$ and the theoretical model is governed by~${\vgx^\vgf}'=\vgM \vgx^\vgf$. Given data~$\vgy^o$, we denote the observation error by~$\egg$, so as~$\egg=\vgy^o-C\vgx^\vgt$. Errors~$\egg$ and~$\ngg$ are null mean white Gaussian noise processes and their respective covariance matrices are~$\vgQ$ and~$\vgR$.

\begin{dfntn}
Given~$\vgx^\vgf(0)$ and $\vgP^\vgf(0)$, the \emph{Kalman-Bucy filter} consists in the two following differential equations, one to estimate the state~$\vgx$ and one for the covariance matrix~$\vgP$, of a differential Riccati type:
\begin{eqnarray*}
\vgx'&=&\vgM\vgx+\vgP C^\text{T}\vgR^{-1}\left(\vgy^o-C\vgx\right),\\
\vgP'&=&\vgM\vgP+\vgP\vgM^\text{T}+ \vgQ-\vgP C^\text{T}\vgR^{-1}C\vgP^\text{T},
\end{eqnarray*}
and the Kalman gain~$\vgK$ is given by~$\vgK=\vgP C^\text{T}\vgR^{-1}.$
\end{dfntn}

Let us explain some links between subsections~\ref{osw} and~\ref{kseek}. In the filters formulations, the feedback is realized thanks to an operator written~$PC^\star TC$ as feedbacks write~$C^\star TC$ and~$P=\text{Id}_{H}$ in subsection~\ref{osw}. Thus one can consider here that~$P$ weights the feedback in comparison to the model. See~\cite{zab76} for some results similar to Proposition~\ref{prop1} about feedbacks~$PC^\star TC$.

We are studying different ways to define the stabilizing operator~$P$, first with the nudging operator~$PC^\star TC=kC^\star C$, where~$k>0$, that is in the framework of Proposition~\ref{prop1}, then with filters. It leads to the following back and forth reconstruction algorithms.

\subsection{Iterative initial data reconstruction methods for inverse problems}\label{iaip}
We get benefits from the reversibility of the wave equations, go back in time and use data again during a backward evolution, from~$t=T$ to~$t=0$, to deduce an approximation of~$f_0$. Such an idea has lead D.~Auroux and J.~Blum to define the Back and Forth Nudging algorithm in~\cite{AB}, then D.~Auroux and E.~Cosme improved it with the Back and Forth SEEK (from private communication, see below for a description of the BF-SEEK). These techniques can be formulated for any sequences of positive operators~$(P^{f}_k),\, (T^{f}_k),\,(P^{b}_k),\, (T^{b}_k)$ as follows:

\begin{dfntn}Given a set of observations~$Cp^o$ and a rough estimate~$p_0$, define~$p^b_0(0)=p_0$. An \emph{iterative reconstruction method} consists in iterating a back and forth process. The forward solution is given by:
\begin{equation}
\label{eqnf}
\partial_{tt}p^f_k=\Delta p^f_k-P^f_k C^\star T^f_k  C(p^f_k-p^o)\text{ in }(0,T)\times\Omega,
\end{equation}
with initial data~$p^f_k(0)=p^f_{k-1}(0)$ and~$p^f_k{'}(0)=0$ and the backward solution solves:
\begin{equation}
\label{eqnb}
\partial_{tt}p^b_k=\Delta p^b_k+P^b_k C^\star T^b_k  C(p^b_k-p^o)\text{ in }(0,T)\times\Omega,
\end{equation}
with \emph{final} data~$p^b_k(T)=p^f_{k}(T)$ and~$p^b_k{'}(T)=p^f_k{'}(T)$. The process is iterated for~$k\geq1$.
\end{dfntn}

Concerning~$PC^\star T C$, it is left implicit that time or space derivatives of~$(p^f_k-p^o)$ and~$(p^b_k-p^o)$ can be considered. The sign preceding~$P^b_k C^\star T^b_k C(p^b_k-p^o)$ changes for convenience of notation as one can notice that, when~$T$ contains a first-order time derivative, then the backward equation~\eqref{eqnb} writes forward:
\begin{equation*}
\label{eqnbf}
\left.
\begin{array}{l}
\partial_{tt}p^b_k=\Delta p^b_k-P^b_k C^\star T^b_k C(p^b_k-p^o),
\end{array}
\right.
\end{equation*}
thanks to the variable substitution~$t\mapsto T-t$. Finally, the correcting term still helps to stabilize the system.

Note that unlike many usual methods, the model is considered here as a weak constraint, which can be useful since it may not be well known (e.g. refer to~\cite{jw} about simulation of inhomogeneous acoustic speed model and relative issues in TAT).

Finally, this kind of algorithms yields successive estimates~$(p^b_k(0))_{k\geq0}$ of the initial object to reconstruct. As explained by Proposition~\ref{prop1} one knows that, under favorable observation conditions, both forward and backward equations are related to exponentially stable semigroups, that leads to the convergence of the algorithm.

\section{Numerical experiments}
\subsection{Discretization and numerical considerations}
Consider the 1-D domain~$\Omega=(-1/2,1/2)$ and~$(0,T)=(0,1)$ both uniformly gridded with steps $\delta x=1/100$ and~$\delta t=1/200$. We deal with the classical finite difference time domain theta-scheme ($\theta$-FDTD) to simulate the wave equation. Sensors are located periodically every~$\delta_{data}$ grid points from the first one.

Upper discussion about observability, stabilizability and controllability have their matching in this situation: a discrete observation condition occur, similar to~\eqref{oc}, which is equivalent too to the discrete stabilizability of the system (we omit details, see~\cite{zua05} and the reference therein when~$\theta=0$). In order to compensate for possibly damaged observation conditions, due either to sensors configuration or to noise, we carry out the solution suggested in~\cite{rtt}, adding an artificial viscous heating attenuation term in the $\theta$-FDTD, that leads to:
$$
\frac{p_{n+1}-2p_n+p_{n-1}}{\delta t^2}=\Delta^\theta_{\delta x}(p_{n-1},p_n,p_{n+1})+\e\Delta_{\delta x}\frac{p_n-p_{n-1}}{\delta t}\pm PC^\star TC\left(\begin{array}{c}p_{n-1}-p^o_{n-1}\\p_n-p^o_n\end{array}\right),
$$
where~$\Delta^\theta_{\delta x}(p_{n-1},p_n,p_{n+1})=\theta\Delta_{\delta x}p_{n-1}+(1-2\theta)\Delta_{\delta x}p_n+\theta\Delta_{\delta x}p_{n+1}$,~$\theta=0.25$,~$\e=0$ or~$\e=(\delta x)^{\alpha}$,~$\alpha\in(1,2)$, and~$\pm$ stands to express both forward and backward implementations.

The term~$\left(\ p_{n-1}-p^o_{n-1}\quad p_n-p^o_n\ \right)^\text{T}$ allows us to consider derivatives of the correcting term in the feedbacks.

\subsection{From Kalman to SEEK filter}
Only main results about the Kalman filter are given and we describe then how to derive the SEEK filter from it (as in~\cite{rbcb}). These algorithms divide into two steps, a \emph{forecast} one and an \emph{analysis} one, in which are taken account the observations to correct the forecast, and two different kinds of parameters are considered, first the states ($\vgx^\vgf$ and $\vgx^\vga$), then the relative error covariance matrices ($\vgP^\vgf$ and $\vgP^\vga$) or their square root ($\vgS^\vgf$ and $\vgS^\vga$). Definitions of the Kalman and SEEK filters are explained respectively on left and right columns below.

\vspace{-1mm}
\begin{dfntn} Given~$\vgx^\vgf_0$,~$\vgP^\vgf_0$ and~$\vgS^\vgf_0=\vgP^{\vgf 1/2}_0$, one iterates:
\vspace{-1mm}
$$
\begin{array}{ccc}
\textbf{Kalman Filter (KF)}&\qquad&\textbf{SEEK Filter}\\
\textit{Analysis step}&\qquad&\textit{Analysis step}\\
\vgK_n=\left[\vgP^{\vgf-1}_n + C^T\vgR^{-1}C\right]^{-1}C^T\vgR^{-1},\hfill&\qquad&\vgK_n=\vgS^\vgf_n\left[I_r+(C\vgS^\vgf_n)^\text{T} \vgR^{-1}(C\vgS^\vgf_n)\right]^{-1} (C\vgS^\vgf_n)^\text{T} \vgR^{-1},\hfill\\
\vgx^\vga_n=\vgx^\vgf-\vgK_n\left[C\vgx^\vgf_n-\vgy^o_n\right],\hfill&\qquad&\vgx^\vga_n=\vgx^\vgf_n - \vgK_n \left[ C\vgx^\vgf_n-\vgy^o_n\right],\hfill\\
\vgP^\vga_n=\left[I_n-\vgK_n C\right]\vgP^\vgf_n,\hfill&\qquad&\vgS^\vga_n=\vgS^\vgf_n \left[ I_r + (C\vgS^\vgf_n)^\text{T} \vgR^{-1} (C\vgS^\vgf_n)\right]^{-1/2},\vspace{1mm}\hfill\\
\textit{Forecast step}&\qquad&\textit{Forecast step}\\
\vgx^\vgf_{n+1}=\vgM\vgx^\vga_{n},\hfill&\qquad&\vgx^\vgf_{n+1}=\vgM \vgx^\vga_i,\hfill\\
\vgP^\vgf_{n+1}=\vgM \vgP^\vga_{n}\vgM^\text{T}+\vgQ.\hfill&\qquad&\vgS^\vgf_{n+1}=\vgM \vgS^\vga_n.\hfill
\end{array}
$$
\noindent where~$I_k$ is the~$k\times k$ identity matrix (in the discrete state space if~$k=n$ and in the reduced rank space if~$k=r$) and~$\vgK_n$ is the Kalman gain, minimizing the trace of the error covariance on~$\vgx^\vga_i$, or the reduced rank Kalman gain in SEEK.
\end{dfntn}

Theoretically, one has~$\vgP^\vgf_{n}=\vgS^\vgf_{n}\vgS^{\vgf T}_{n}+\vgQ$ in SEEK, but since~$\vgQ$ is not well known, we define~$\vgP^\vgf_{n}=(1+\gamma)\vgS^\vgf_{n}\vgS^{\vgf T}_{n}$,
where~$\gamma\in(0,1)$, which avoids an additional decomposition in SEEK. It is set similarly in KF.

In KF, if the state space has a range of~$n$, the forecast step necessitates~$2n$ model evolution steps to get the forecast error covariance matrix from the analysis error covariance matrix, such that a less optimal gain may be considered to reduce calculation cost. This is the purpose of the SEEK to yield such a gain by considering a reduced rank gain more simple to obtain. Since the error covariance matrices are symmetric and positive-definite, Pham et al.~\cite{pham} suggested to consider the following decomposition: if~$P$ is a symmetric positive-definite~$n\times n$ matrix whose~$r\leq n$ larger eigenvalues are~$\lambda_1\geq\ldots\geq \lambda_r>0$ and their corresponding eigenvectors~$V_1,\ldots,V_r$, the \emph{reduced decomposition} of~$P$ is defined as the~$n\times r$ matrix~$S=[\sqrt{\lambda_1}V_1 \ldots \sqrt{\lambda_r}V_r]$. Since SEEK consists in using the reduced decomposition~$S^f$ of~$P^f$,~$r$ model runs define the forecast error covariance matrix.

\subsection{Results}
Four methods are compared: Time Reversal (TR), Back and Forth Nudging (BFN), Back and Forth SEEK (BF-SEEK) and Kalman Filter (KF). Let us set~$\vgR=0.3I$,~$\gamma=0.01$ for both BF-SEEK and KF, and the reduced rank~$r=120$, which is reduced when~$\text{Rank}\left[ I + (C\vgS^\vgf_i)^\text{T} \vgR^{-1} (C\vgS^\vgf_i)\right]<r$. An additive white Gaussian noise of level~$\nu$ is added to the data when said. The nudging gain is set to~$0.9\delta t$ for implementation limitations.

\begin{figure}[!h]
\centering
\includegraphics[width=.8\textwidth,height=.13\textwidth]{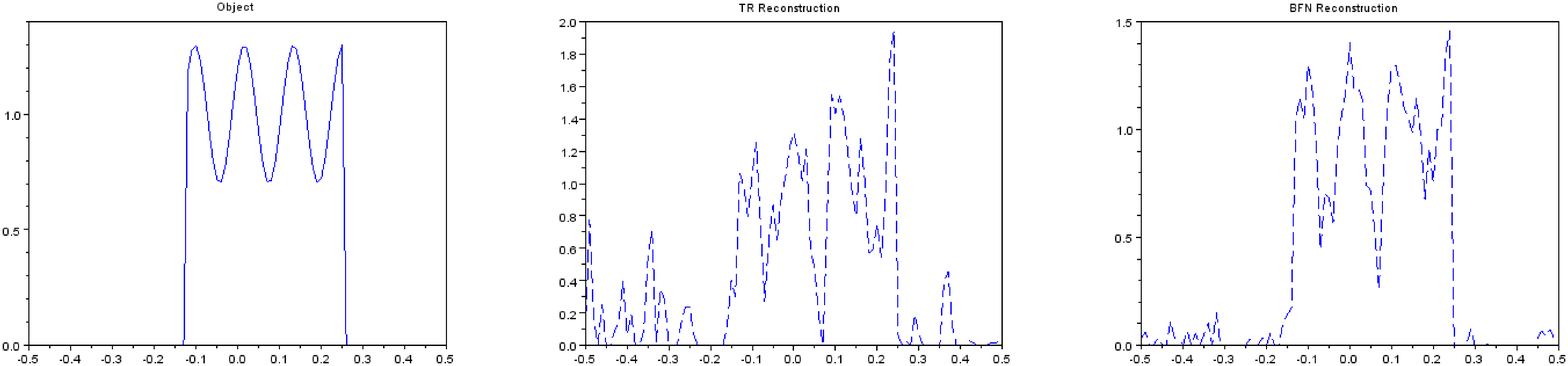}\\
\includegraphics[width=.53\textwidth,height=.13\textwidth]{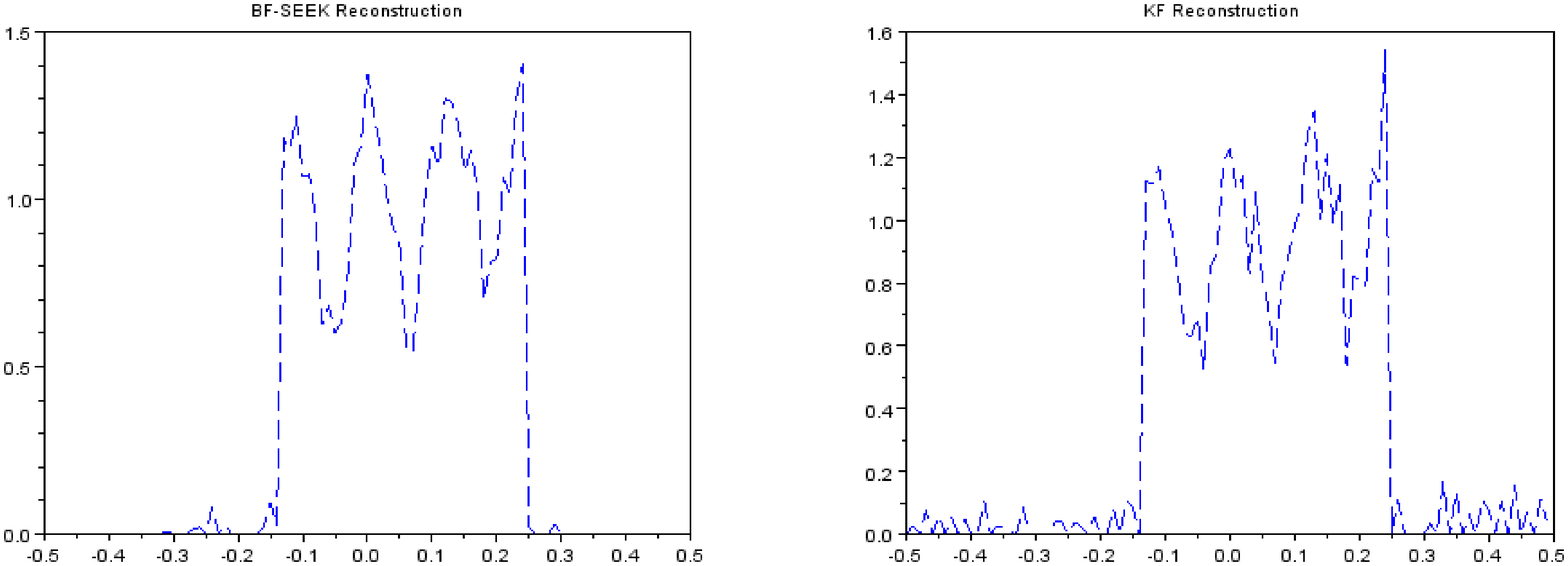}
\caption{\label{fig2}Object to reconstruct and reconstructions by TR, BFN, BF-SEEK and KF for~$\delta_{data}=10$ and~$\nu=30\%$.}
\end{figure}

\begin{figure}[!h]
\centering
\includegraphics[width=.2\textwidth,height=.15\textwidth]{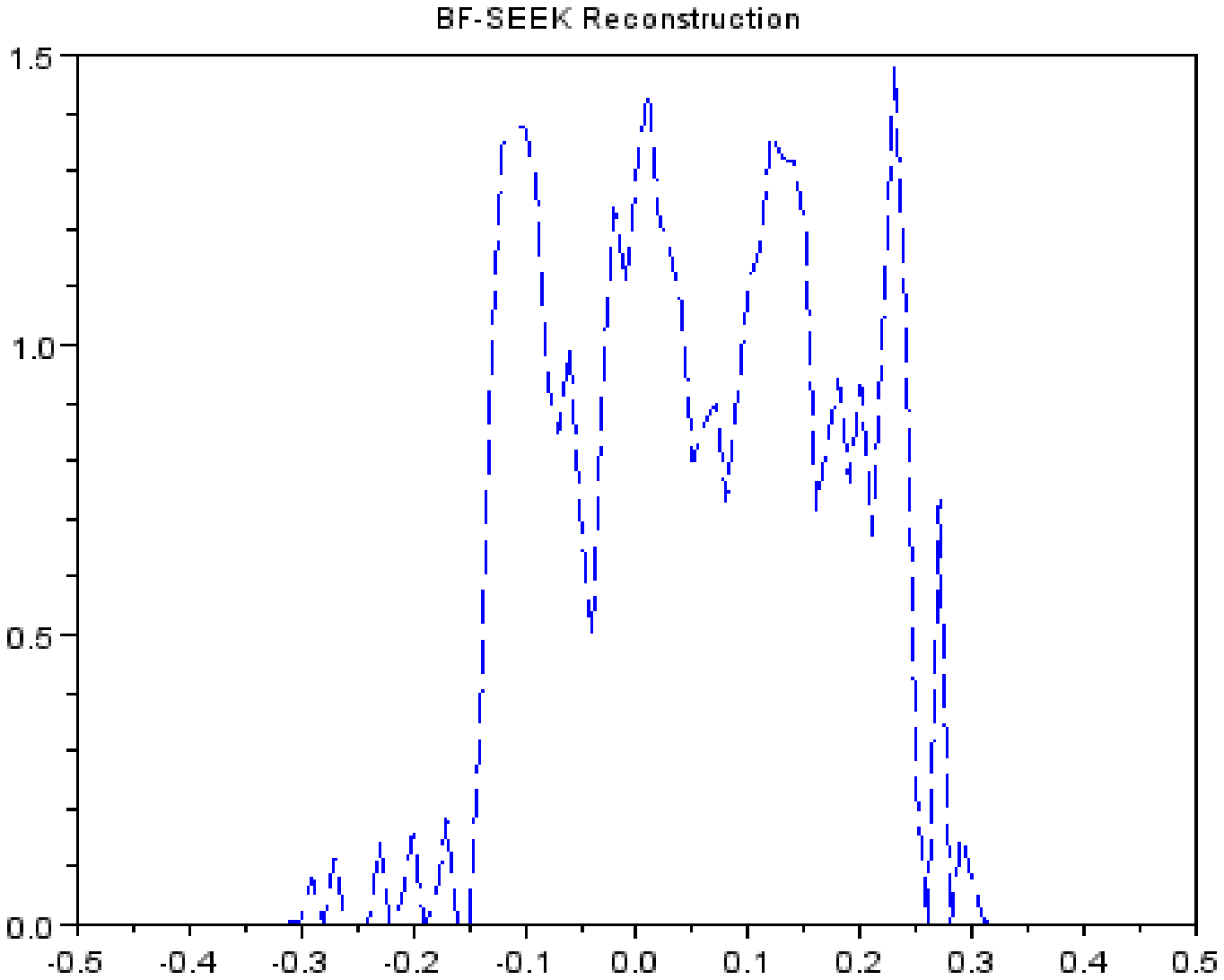}
\includegraphics[width=.2\textwidth,height=.147\textwidth]{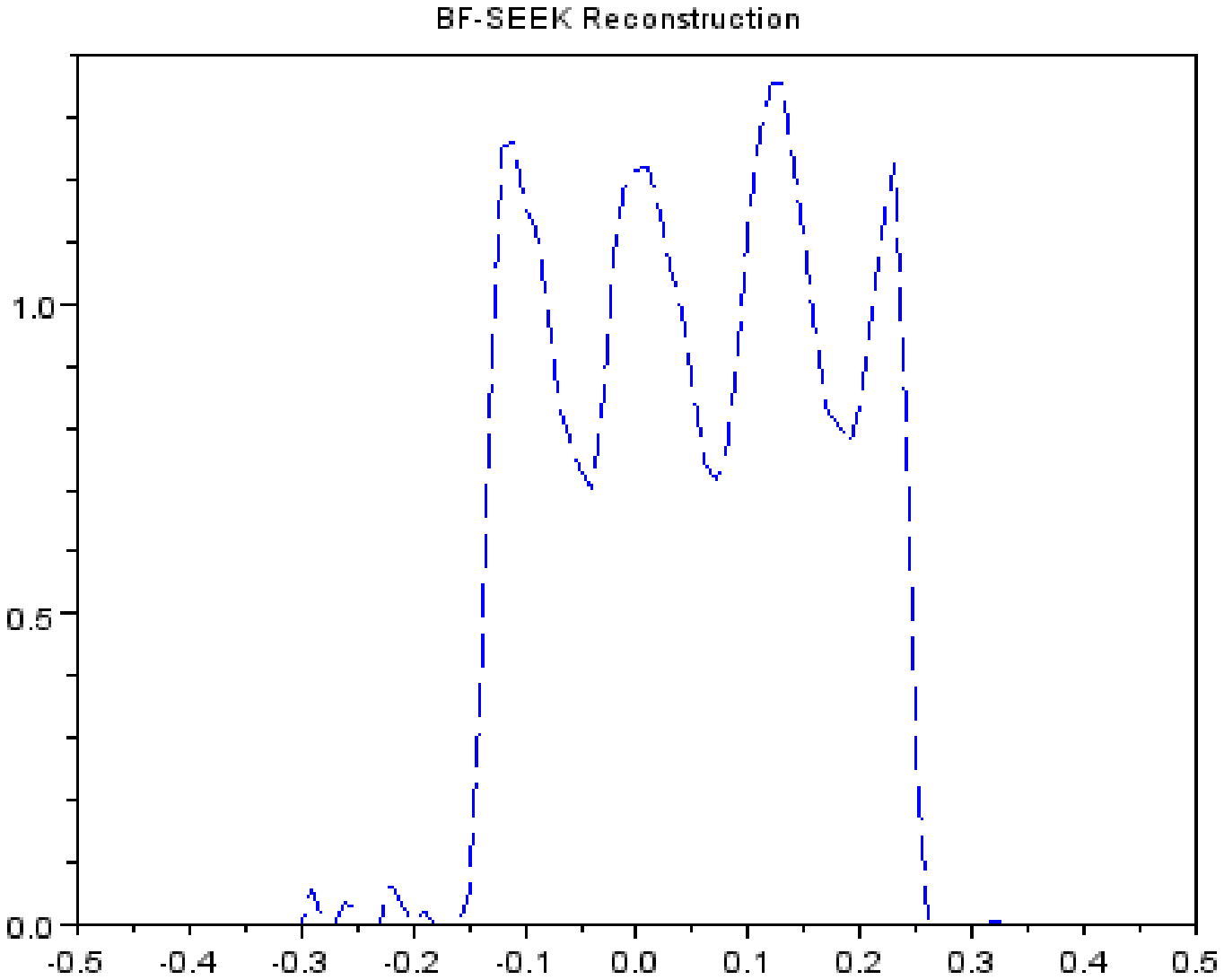}
\includegraphics[width=.2\textwidth,height=.15\textwidth]{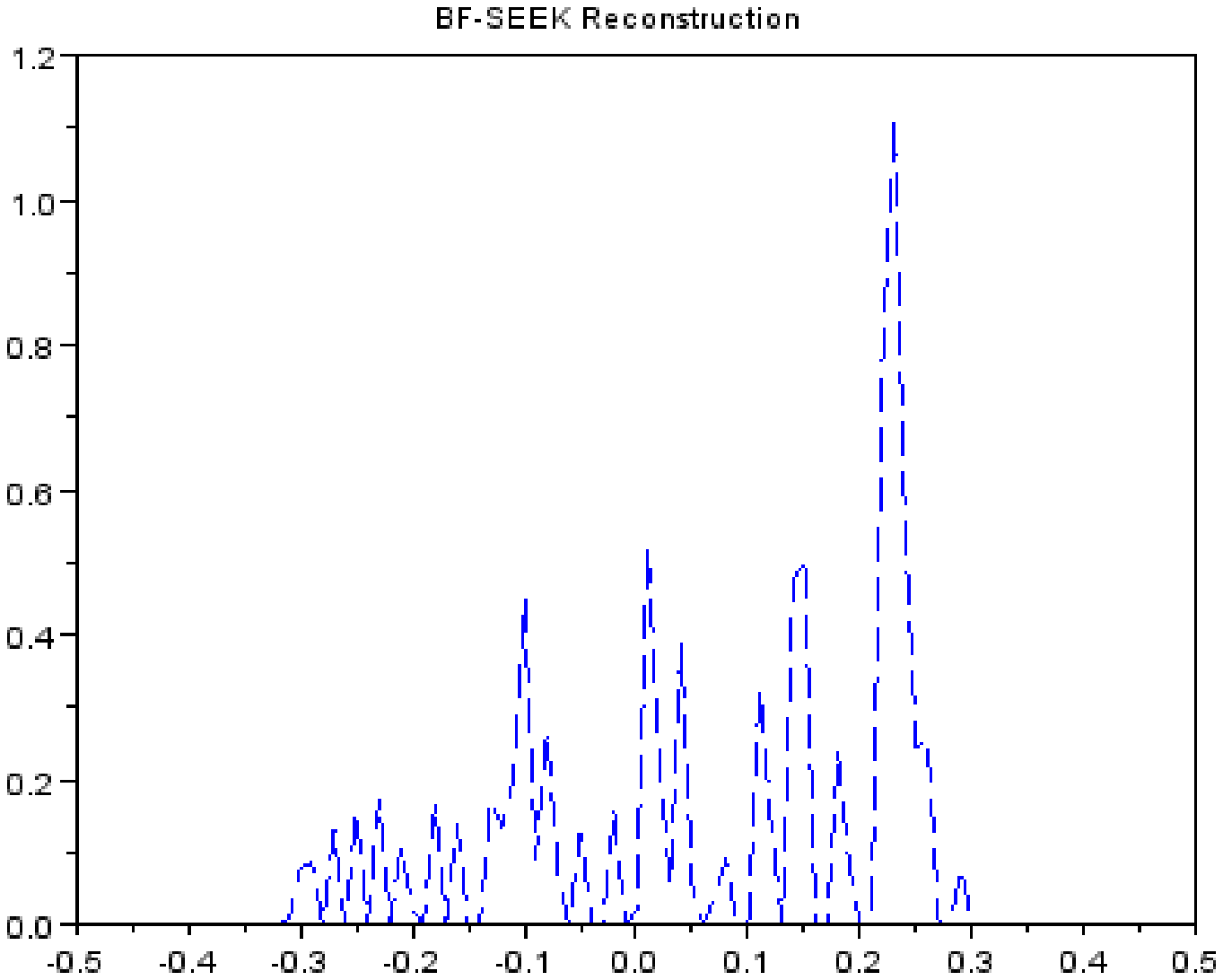}
\includegraphics[width=.2\textwidth,height=.147\textwidth]{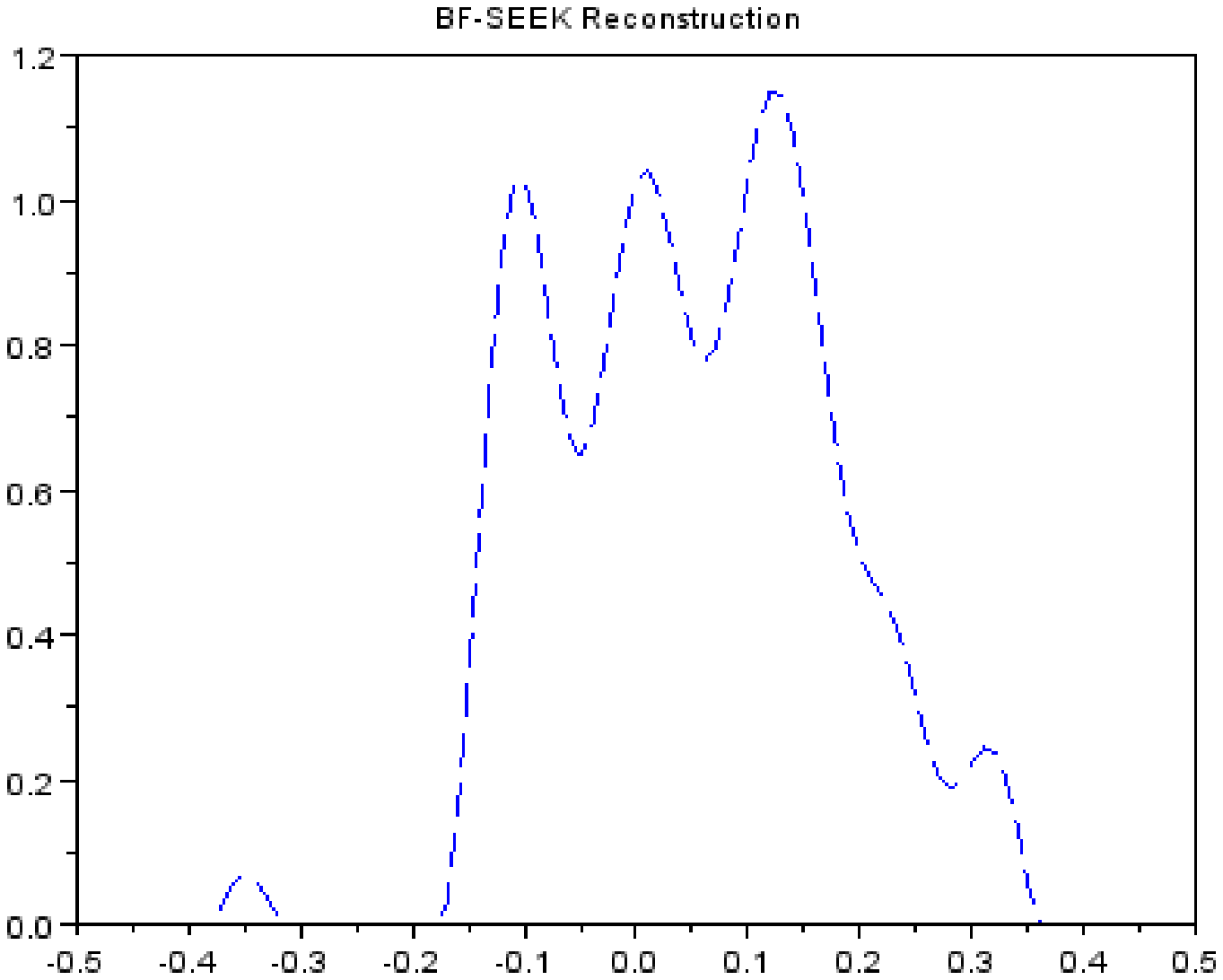}
\caption{\label{fig3}BF-SEEK reconstructions. Case~$\delta_{data}=99$ and $\nu=30\%$: $\alpha=0$~(a) and~$\alpha=2$~(b). Case~$\delta_{data}>99$: $\alpha=0$~(c) and~$\alpha=1.8$~(d).}
\end{figure}

Table~\ref{fig1} and Fig.~\ref{fig2} and~\ref{fig3} show RMS errors and some relative reconstructions obtained in various situations. The object to reconstruct is shown on the upper left part of figure~\ref{fig2} then, to the right, one sees TR and BFN reconstructions, followed by BF-SEEK and KF reconstructions.
\vspace{-.4cm}
\begin{table}[h!]
\centering
$$
\begin{array}{lcccc}
\centering\text{Settings}&\text{TR}&\text{BFN}&\text{BF-SEEK}&\text{KF}\\
\delta_{data}=10&5.4&0.5\; (100)&0.6\; (3)&6\\
\delta_{data}=10,\ \nu=30\%&38.8&16.3\; (5)&10.3\; (2)&15.7\\
\delta_{data}=99&8.3&2.6\; (100)&6.3\; (1)&44.2\\
\delta_{data}=99,\ \nu=30\%&14.1&10.8\; (14)&19.9\; (1)&46.4\\
\delta_{data}=99,\ \nu=30\%,\ \alpha=2&9.8&20.7\; (8)&12.2\; (1)&42\\
\delta_{data}>99&51.1&20.7 (7)&\begin{array}{c} 83.3\; (3)\\31.8\; \text{if}~\alpha=1.8\;(3)\end{array}&95
\end{array}
$$
\caption{\label{fig1} RMS errors and number of iterations to converge (in brackets) for~10 ($\delta_{data}=10$),~2 ($\delta_{data}=99$) and~1 ($\delta_{data}>99$) sensor(s), noisy ($\nu=30\%$) or noiseless data, with or without the numerical attenuation $\delta x^\alpha$.}
\vspace{-.5cm}
\end{table}


As KF reacts quite well to noise addition, it shows much more sensitivity to the number of sensors and easily fails. BF-SEEK offers obvious improvements for both calculation cost and reconstruction error (Fig.~\ref{fig2}). When~2 sensors are left, one can observe possible interesting effects of the attenuation term against noise with TR and BF-SEEK. Nevertheless, it may damage the reconstruction (with BFN) or have insignificant consequences (with TR). When only~1 sensor is left, only BFN and BF-SEEK are robust enough to yield a good approximation of the object, but BF-SEEK needs to be corrected with the attenuation to keep stable (Fig.~\ref{fig3}).

\section{Conclusions}
A common formulation for iterative stabilization of reversible evolution systems is given. It is used to define methods which solve some inverse problems for wave equations. Experiments show that the techniques we introduced may offer an alternative to usual inverse methods for TAT problem. In applications, knowledge about the quality of the sensors, the background and the model approximations can be used by filters, but in this case they need to be precisely tuned.

The use of an artificial attenuation is motivated by good results, and allows one to consider lossy medium in back and forth implementations when the physical loss does not exceed the numerical attenuation.

When space dimension increases, we first face an excessive calculation cost (one BF-SEEK iteration with rank~60 equals almost one thousand BFN in computation time). So one would get interested in hybrid methods, e.g. getting first estimates with TR and BFN and then reconstructing the solution with filters. A solution is offered by Perfectly Matched Layers to reduce the space domain of calculation and for which a first order discrete scheme formulation is necessary for the filters.


\end{document}